\documentstyle[a4,leqno,amsfonts,12pt]{article}
\newcommand{\C}{{\bf C}}

\newcommand{\OC}{\overline{{\bf C}}}

\newcommand{\N}{{\bf N}}

\newcommand{\bP}{{\bf P}}
\newcommand{\D}{{\bf D}}

\newcommand{\de}{\delta}

\newcommand{\mb}{\mbox}

\newcommand{\beq}{\begin{equation}}
\newcommand{\eeq}{\end{equation}}
\newcommand{\oge}{\succeq}
\newcommand{\ole}{\preceq}
\newcommand{\ve}{\varepsilon}

\newcommand{\ov}{\overline}
\newcommand{\al}{\alpha}
\newcommand{\be}{\beta}

\newcommand{\Om}{\Omega}

\newcommand{\z}{\zeta}

\newcommand{\la}{\lambda}

\newtheorem{th}{Theorem}

\newtheorem{lem}{Lemma}

\newcommand{\ueberschrift}{\bigskip\goodbreak\noindent\bigskip}
\newcounter{theabsatz}
\newcommand{\absatz}[1]{\stepcounter{theabsatz} \ueberschrift
               {\large \bf \arabic{theabsatz}. {#1}} \setcounter{equation}{0}}

\begin{document}

\begin{center}
{\bf \large ON THE CHRISTOFFEL FUNCTION FOR THE
GENERALIZED JACOBI MEASURES ON A QUASIDISK}\\[2ex]

 Vladimir V.
Andrievskii \\[3ex] 

\end{center}

\begin{abstract}
We establish the exact (up to the constants) double inequality for the Christoffel function
for a measure supported on a Jordan domain bounded by a quasiconformal
curve. We show that  this quasiconformality of the boundary 
 cannot be omitted.
\end{abstract}

\footnotetext{
Date received:$\quad$. Communicated by


 }

\absatz{Introduction and Main Results}

Denote by $\bP_n$ the set of all complex polynomials of degree at most $n\in\N_0:=\{0,1,2,\ldots\}$.
For a
 finite Borel measure $\nu$ on the complex plane $\C$ such that its support is compact
and it consists of infinitely many points,
 and a parameter $1\le p<\infty$,
the {\it $n$-th Christoffel function} associated with $\nu$ and $p$, is defined by 
\beq\label{1.1}
\la_n(\nu,p,z):=\inf _{p_n\in \bP_n \atop p_n(z)=1}\int|p_n|^pd\nu,\quad z\in\C.
\eeq
This function plays an important role in the theory of orthogonal polynomials,
in particular, due to the following {\it Christoffel Variational Principle}
(see \cite[p. 78]{totsur} or \cite[p. 309]{sim}): 
\beq\label{1.2}
\la_n(\nu,2,z)=\left(\sum_{j=0}^n|\pi_j(\nu,z)|^2\right)^{-1}, \quad 
z\in\C,
\eeq
where $\pi_j(\nu,\cdot)$ is the $j$-th orthonormal polynomial
associated with measure $\nu$.

The starting point of our consideration consists of two groups of results.
The first group includes recent findings in \cite{tot10}-\cite{var} about the behavior 
of $\la_n(\nu,p,z)$ in the case where $\nu$ is supported on a Jordan arc or curve.
The second group includes results  in \cite{sue, abd, abddeg, gupusast, sty} 
 about the behavior of $\pi_n(\nu,z)$ in the case of a (weighted) area type measure $\nu$.
We refer the reader to these papers 
for the
further references.

We consider measures supported on the closure $\ov{G}$
of a domain $G\subset\C$ bounded by a Jordan curve $L:=\partial
G$. Let
$\Om:=\OC\setminus\ov{G}$, where $\OC:=\C\cup\{\infty\}$ is the extended complex plane.
The  Riemann mapping function $\Phi:\Om\to\D^*:=\{w:|w|>1\}$ normalized by
$$
\Phi(\infty)=\infty,\quad \Phi'(\infty):=\lim_{z\to \infty}\frac{\Phi(z)}{z}>0
$$
plays an essential role in our consideration, which from this point of view, can be compared with
the results in the above mentioned papers.

We focus our attention  to the case where $G$ is a bounded {\it quasidisk}, i.e., $L$
is a bounded {\it qusiconformal curve} (see \cite{ahl, lehvir}) which geometrically means
that for every pair of different points $z_1,z_2\in L$,
\beq\label{z1}
\min(\mb{diam }L', \mb{diam }L'')\le C_L|z_1-z_2|,
\eeq
where $L'$ and $L''$ denote the two connected components (subarcs)
of $L\setminus\{z_1,z_2\}$,
 diam $S$ is the diameter of a set $S\subset\C$, and $C_L\ge1$ is a constant depending only on $L$.

For fixed $z_j\in L$ and $\al_j>-2, j=1,\ldots,m$, consider the  {\it weight function}
\beq\label{1.3}
h(z):=
\left\{\begin{array}{ll}
   \displaystyle
 h_0(z)\prod_{j=1}^m|z-z_j|^{\al_j}
 &\mb{ if } z\in G,
\\[2ex]
0&\mb{ if }z\in\C\setminus G,
\end{array}\right.
\eeq
where, for a  measurable function $h_0$, the inequality
$$
C_h^{-1}\le h_0(z)\le C_h,\quad z\in G
$$
holds with a constant $C_h\ge1$ depending only on $h$.

A measure $\nu$ supported on $\ov{G}$ and determined by $d\nu=hdm$, where $dm$
stands for the $2$-dimensional Lebesgue measure (area) in the plane,
 is called the
{\it generalized Jacobi measure}. 

Let 
$$
d(z,S):=\mb{dist}(\{z\},S):=\inf_{\z\in S}|z-\z|,\quad
z\in\C, S\subset\C,
$$
and let for $\de>0$ and $z\in L$,
$$
L_\de:=\{\z\in\Om:|\Phi(\z)|=1+\de\},\quad \rho_{\de}(z):= d(z,L_\de).
$$
\begin{th}\label{th1}
Let $G$ be a quasidisk, $\nu$ be the generalized Jacobi measure, and let
$1\le p<\infty$. Then for $n\in\N:=\{1,2,\ldots\}$ and $z\in L$,
\beq\label{1.6}
C^{-1}\le\la_n(\nu,p,z)\rho_{1/n}(z)^{-2}\prod_{j=1}^m
(|z-z_j|+\rho_{1/n}(z))^{-\al_j}\le C
\eeq
holds with $C=C(G,h,p)\ge1$.
\end{th}
According to (\ref{1.2}) and (\ref{1.6}), for the orthogonal polynomials $\pi_n(\nu,z)$
and $z\in L$, we have
\begin{eqnarray}
|\pi_n(\nu,z)|&\le &\la_n(\nu,2,z)^{-1/2}
\nonumber\\
\label{1.7}
&\le&
 C^{1/2}
\rho_{1/n}(z)^{-1}\prod_{j=1}^m
(|z-z_j|+\rho_{1/n}(z))^{-\al_j/2}.
\end{eqnarray}
Using (\ref{1.7}) and well-known distortion properties of  conformal mappings
with quasiconformal extension (see \cite{pom75}) one can obtain  more specialized
bounds for orthogonal polynomials which can be found, for example, in
\cite{sue, abd, abddeg} where they are proved by other methods.

If $G=\D:=\{z:|z|<1\}$ and $d\nu = h_0dm$, then (\ref{1.7}) becomes
\beq\label{1.nn}
|\pi_n(\nu,z)|\le C^{1/2} n,\quad z\in\ov{\D}.
\eeq
Keeping in mind the Rakhmanov's \cite{rak} solution of the Steklov problem 
(for more details, see  \cite{aptdentul} or \cite{den}), it is tempting to conjecture
that (\ref{1.nn}) as well as (\ref{1.7}) cannot be improved.

The inequality (\ref{1.6}) can also be used to estimate $|\pi_n(\nu,z)|$
from below. For example,
if $\al_j=0$, i.e., $h(z)=h_0(z)$  for all $z\in G$,
then, by virtue of (\ref{1.6}),
for any quasidisk $G$ and $1\le p<\infty$ we have
\beq\label{1.8}
C^{-1}\rho_{1/n}(z)^2\le \la_n(\nu,p,z)\le C\rho_{1/n}(z)^2,\quad
z\in L,n\in\N,
\eeq
which, together with (\ref{2.4z}) below, imply
that there exists $k=k(G)\in\N\setminus\{1\}$ 
such that for $z\in L$ and $ n\in\N$,
$$
\sum_{j=n+1}^{kn}|\pi_j(\nu,z)|^2=
\la_{kn}(\nu,2,z)^{-1}-\la_{n}(\nu,2,z)^{-1}
\ge \frac{1}{2C\rho_{1/(kn)}(z)^2}.
$$
Therefore,
$$
\max_{n<j\le kn}|\pi_j(\nu,z)|\ge\frac{\ve}{\sqrt{n}\rho_{1/(kn)}(z)},
\quad \ve:=(2kC)^{-1/2},
$$
that is,
$$
\max_{n<j\le kn}(\sqrt{j} \rho_{1/j}(z)|\pi_j(\nu,z)|)\ge\ve,
$$
which yields that for any $z\in L$ there exists an infinite set $\Lambda_z\subset\N$ such that
\beq\label{*}
|\pi_n(\nu,z)|\ge\frac{\ve}{\sqrt{n}\rho_{1/n}(z)},\quad n\in\Lambda_z.
\eeq
Note that the case $h(z)\equiv 1$ on $G=\D$ shows the exactness of
(\ref{*}) (up to the constant $\ve$).

Next, consider the domain
$$
G^*:=\{z=x+iy:0<x<1, |y|<e^{-1/x}\}
$$
which is obviously not a quasidisk.
\begin{th}\label{th2}
For the area  measure $m^*$  supported on $\ov{G^*}$, $1\le p<\infty$, and 
$k\in\N$,
\beq\label{1.9}
\lim_{n\to\infty}\frac{\la_n(m^*,p,0)}{\rho_{1/n}(0)^k}=0.
\eeq
\end{th}
Comparing the left-hand side of
(\ref{1.8}) and  (\ref{1.9}) shows that the requirement in Theorem \ref{th1} on $G$
to be a quasidisk cannot be omitted.

 Using the approach from the proof of \cite[Theorem 2]{and05}
or \cite[Corollary 2.5]{var} the same inequality
(\ref{1.6}) can be proved if $G$ is replaced by a finite union
of quasidisks lying exterior to each other. We do not dwell on this purely technical
problem.

The structure of this paper is as follows. Section 2 contains auxiliary results
from theory of quasiconformal mappings and constructive function theory in the 
complex plane. In Section 3, we prove the main results, i.e., Theorem \ref{th1}
and Theorem \ref{th2}.

In what follows, we always assume that $G$ is a quasidisk and $h$ is a
generalized Jacobi measure.
We  use
the convention that $c,c_1,\ldots
$ denote positive constants and
$\ve,\ve_1,\ldots$ sufficiently small positive constants
(different in different sections). If not stated otherwise, we
assume that these constants can depend only on $G,p$, and $h$.
 For the nonnegative functions $a$
and $b$ we write $a\preceq b$ if $a\le cb$, and $a\asymp b$ if
$a\preceq b$ and $b\preceq a$ simultaneously.

We complete this section with the additional notation:
$$
D(z,r):=\{\z:|\z-z|<r\},\,\,C(z,r):=\{\z:|\z-z|=r\},\quad z\in\C,r>0.
$$

\absatz{Auxiliary Results and Constructions}

We begin with estimation of two integrals. 
\begin{lem}\label{lem2.1}
Let $\de>0,\al>-2,\be>2+|\al|$. Then for $z',z''\in\C$
we have
\beq\label{2.1}
I:=
\int_{\C}(|\z-z''|+\de)^{-\be}|\z-z'|^\al dm(\z)
\le c_1\de^{2-\be}(|z'-z''|+\de)^\al,
\eeq
where $c_1=c_1(\al,\be)$.
\end{lem}
{\bf Proof.}
Consider two particular cases.

If $|z'-z''|\le\de$, then using the polar coordinates with center at $z'$,
we obtain
\begin{eqnarray}
I&\le&\de^{-\be}\int_{D(z',2\de)}|\z-z'|^\al dm(\z)
+ \int_{\C\setminus D(z',2\de)}|\z-z'|^{\al-\be} dm(\z) 
\nonumber\\
\label{2.2}
&\le& 2\pi\de^{-\be}\int_0^{2\de}r^{\al+1}dr
+2\pi\int_{2\de}^\infty r^{\al-\be+1}dr
\ole \de^{\al-\be+2}.
\end{eqnarray}
If $|z'-z''|>\de$, then letting $d:=|z'-z''|$,
$$
D':=D\left(z',\frac{d}{2}\right), \quad D'':=D\left(z'',\frac{d}{2}\right),
$$
$$
U_1:=D(z',2d)\setminus(D'\cup D''),\quad
U_2:=\C\setminus D(z',2d),
$$
and using the polar coordinates with centers at $z'$ and
$z''$ respectively, we have
\begin{eqnarray}
I&\ole&d^{-\be}\int_{D'}|\z-z'|^\al dm(\z)
+d^{\al}\int_{D''}(|\z-z''|+\de)^{-\be} dm(\z)
\nonumber\\
&&+ d^{\al-\be}\int_{U_1}dm(\z)+\int_{U_2}|\z-z'|^{\al-\be}dm(\z)
\nonumber\\
&\le& 2\pi d^{-\be}\int_0^{d/2}r^{\al+1}dr+2\pi d^\al\int_{0}^{d/2}(r+\de)^{-\be+1}dr
\nonumber\\
&&+d^{\al-\be} \pi 4 d^2+ 2\pi\int_{2d}^\infty r^{\al-\be+1}dr
\nonumber\\
\label{2.3}
&\ole& d^{-\be+\al+2}+d^\al\de^{2-\be}
\asymp d^\al\de^{2-\be}.
\end{eqnarray}
Comparing (\ref{2.2}) and (\ref{2.3})  we obtain (\ref{2.1})

\hfill$\Box$

\begin{lem}\label{lem2.2}
Let $0<\de<\ve$, $\al_j>-2,j=1,\ldots,m$,
and $\be>2+\sum_{j=1}^m|\al_j|$. Suppose that
points
$z_1,\ldots, z_m\in \C$ satisfy
$$
|z_j|<c,\quad |z_j-z_k|>4\ve,\quad j\neq k.
$$
Then, for any $z\in\C$ with $|z|<c$, we have
\begin{eqnarray}
I^*(z)&:=&\int_{D(0,c)}(|\z-z|+\de)^{-\be}\prod_{j=1}^m|\z-z_j|^{\al_j} dm(\z)
\nonumber\\
\label{2.4}
&\le&
c_2\de^{2-\be}\prod_{j=1}^m(|z-z_j|+\de)^{\al_j},
\end{eqnarray}
where $c_2=c_2(\al_1,\ldots,\al_m,z_1,\ldots,z_m,\ve, c, \be).$
\end{lem}
{\bf Proof}. Let $\al:=\sum_{j=1}^m|\al_j|$ and
$$
D_j:=D(z_j,2\ve),D_j':=D(z_j,\ve),
\quad j=1,\ldots,m.
$$
Consider two particular cases.

If $z\not\in\cup_{j=1}^m D_j$, then
using the polar coordinates with centers at $z_j$ and $z$ respectively
we have
\begin{eqnarray}
I^*(z)&\ole& \sum_{j=1}^m\int_{D_j'}|\z-z_j|^{\al_j}dm(\z)+
\int_\C(|\z-z|+\de)^{-\be}dm(\z)\nonumber\\
\label{2.5}
&\ole& \sum_{j=1}^m\int_0^\ve r^{\al_j+1}dr
+\int_0^\infty(r+\de)^{-\be+1}dr\ole\de^{2-\be}.
\end{eqnarray}
If $z\in D_k$ for some $k=1,\ldots,m$, then, 
applying Lemma \ref{lem2.1}, we obtain
\begin{eqnarray}
I^*(z)&\ole&\sum_{j=1\atop j\neq k}^m\int_{D_j'}|\z-z_j|^{\al_j}dm(\z)
+\int_{\C}(|\z-z|+\de)^{-\be}|\z-z_k|^{\al_k}dm(\z)
\nonumber\\
\label{2.61}
&\ole&
\de^{2-\be}\prod_{j=1}^m(|z-z_j|+\de)^{\al_j}.
\end{eqnarray}
Comparing (\ref{2.5}) and (\ref{2.61})  we have (\ref{2.4}).

\hfill$\Box$

Now let $z\in L=\partial G, 0<r\le\de<($diam $G)/4$, and $\al>-2$.
Since by the definition of a quasiconformal curve (\ref{z1})
$$
|G\cap C(z,r)|\oge r,
$$
where $|S|$ means the {\it linear measure}, i.e. {\it length}, of $S\subset\C$,
we have 
\beq\label{2.6}
\int_{G\cap D(z,\de)}|\z-z|^\al dm(\z)=\int_0^\de r^\al
|G\cap C(z,r)| dr\oge \de^{\al+2}.
\eeq
Therefore, if $Z:=\{z_1,\ldots,z_m\}\subset L$ and
$\al_j>-2,j=1,\ldots,m$ are fixed, then for $z\in L$ and
$\de<\min_{j\neq k}|z_j-z_k|/4$,
\beq\label{2.7}
\int_{G\cap D(z,\de)}\prod_{j=1}^m|\z-z_j|^{\al_j}dm(\z)
\ge \ve_1\de^2\prod_{j=1}^m(|z-z_j|+\de)^{\al_j},
\eeq
where $\ve_1=\ve_1(G,Z,\al_1,\ldots,\al_m)$.

Indeed, let $d:=d(z,Z)=|z-z_k|$ for some $k=1,\ldots,m$.
If $\de\ge 2d$, then by virtue of (\ref{2.6})
\begin{eqnarray}
A&:=&
\int_{G\cap D(z,\de)}\prod_{j=1}^m|\z-z_j|^{\al_j}dm(\z)
\oge\int_{G\cap D(z_k,\de/2)}|\z-z_k|^{\al_k}dm(\z)
\nonumber\\
\label{2.8}
&\oge& \de^{\al_k+2}\asymp\de^2\prod_{j=1}^m
(\de+|z-z_j|)^{\al_j}=:B.
\end{eqnarray}
If $\de<2d$, then, according to (\ref{2.6}),
\beq\label{2.9}
A\oge d^{\al_k}\int_{G\cap D(z,\de/4)}dm(\z)
\oge\de^2 d^{\al_k}\asymp B.
\eeq
Comparing (\ref{2.8}) and (\ref{2.9})  we have (\ref{2.7}).

Next, we introduce  auxiliary families of quasiconformal curves and mappings
as follows. 
Let $K\ge 1$ be a coefficient  of quasiconformality of  $L$. It is well known 
(see \cite[Chapter IV]{ahl}) that the Riemann mapping function $\Phi$
can be extended to a $K^2$-quasiconformal homeomorphism $\Phi:\OC\to\OC$. Hence,
each curve
$$
L^*_\de:=\{z:|\Phi(z)|=1-\de\},\quad 0\le\de<1
$$
is $K^2$-quasiconformal. Denote by $\Om^*_\de$ the unbounded connected component
of $\OC\setminus L^*_\de$. The Riemann conformal mapping $\Phi_\de:\Om^*_\de\to
\D^*$ with the normalization
$$
\Phi_\de(\infty)=\infty,\quad \Phi_\de'(\infty)>0
$$
can be extended to a $K^4$-quasiconformal homeomorphism 
$\Phi_\de:\OC\to\OC$. 
Note that $\Phi_0=\Phi,\Psi_0=\Psi,$ and $L_0^*=L$.
To study  metric properties of $\Phi_\de$
and $\Psi_\de:=\Phi_\de^{-1}$, we use the following statement.
\begin{lem}\label{lem2.3}
(see \cite[p. 97, Theorem 4.1]{andbeldzj} or \cite[p. 29, Theorem 2.7]{andbla})
Suppose that $F:\OC\to\OC$ is a $Q$-quasiconformal mapping with 
$Q\ge 1$ and $F(\infty)=\infty$.
Assume also that $\z_j\in\C,w_j:=F(\z_j),j=1,2,3$. Then:

(i) the conditions $|\z_1-\z_2|\le c_3|\z_1-\z_3|$ and
$|w_1-w_2|\le c_4|w_1-w_3|$ are equivalent; besides, the constants $c_3$  and $c_4$
are mutually dependent and dependent on $Q$;

(ii) if $|\z_1-\z_2|\le c_3|\z_1-\z_3|$, then
$$
c_5^{-1}\left|\frac{w_1-w_3}{w_1-w_2}\right|^{1/Q} \le
 \left|\frac{\z_1-\z_3}{\z_1-\z_2}\right|\le
c_5\left|\frac{w_1-w_3}{w_1-w_2}\right|^Q,
$$
where $c_5=c_5(c_3,Q).$
\end{lem}
Let
$$
\tilde{\z}_\de:=\Psi((1+\de)\Phi(\z)),\quad \z\in\ov{\Om}\setminus\{\infty\}, \de>0.
$$
For $z\in L$ and $\de>0$, let a point $z^*_\de\in L_\de$ satisfy $|z-z^*_\de|
 =\rho_\de(z)$. Applying Lemma \ref{lem2.3} with $F=\Phi$ and the triplet of points
 $z,\tilde{z}_\de, z_\de^*$ we have
\beq\label{2.2z} 
\rho_\de(z)\asymp|z-\tilde{z}_\de|,\quad z\in
L.
\eeq
Moreover, we claim that for $z,\z\in L$
and $0<\de\le 1$,
\beq\label{2.nn}
|z-\tilde{\z}_\de|\asymp |z-\z|+\rho_\de(z).
\eeq
Indeed, by Lemma \ref{lem2.3} with $F=\Phi$,
$$
|z-\z|\ole|z-\tilde{\z}_\de|\quad \mb{and}\quad |z-\tilde{z}_\de|\ole |z-\tilde{\z}_\de|,
$$
i.e.,
$$
|z-\z|+\rho_\de(z)\asymp |z-\z| + |z-\tilde{z}_\de|\ole |z-\tilde{\z}_\de|.
$$
Furthermore, the same Lemma \ref{lem2.3} with $F=\Phi$ also implies that
if $|\Phi(z)-\Phi(\z)|>\de$ then
$$
|z-\tilde{\z}_\de|\asymp|z-\z|\oge |z-\tilde{z}_\de|
$$
as well as if $|\Phi(z)-\Phi(\z)|\le \de$ then
$$
|z-\tilde{\z}_\de|\asymp|z-\tilde{z}_\de|\oge |z-\z|.
$$
That is, in both cases we have
$$
|z-\tilde{\z}_\de|\oge |z-\z|+|z-\tilde{z}_\de|
\asymp |z-\z|+\rho_\de(z)
$$
which completes the proof of (\ref{2.nn}).

 Next,
for $0<v<u\le 1$ and $z\in L$,
 Lemma \ref{lem2.3} with $F=\Phi$ and the triplet of points
 $z,\tilde{z}_v,\tilde{z}_u$ as well as (\ref{2.2z}) imply
 \beq\label{2.4z}
\left(\frac{u}{v}\right)^{1/K^2}\ole\frac{\rho_u(z)}{\rho_v(z)}\ole
\left(\frac{u}{v}\right)^{K^2}.
\eeq
For $\xi\in\Om\setminus\{\infty\}$, let
$
\xi_L:=\Psi(\Phi(\xi)/|\Phi(\xi)|)
$
and let $\xi_*\in L$ satisfy $d(\xi,L)=|\xi-\xi_*|.$
Applying Lemma \ref{lem2.3} with $F=\Phi$ and the triplet 
of points $\xi, \xi_L,\xi_*$
we obtain
\beq\label{2.5z}
|\xi-\xi_L|\asymp d(\xi,L).
\eeq
Therefore, for $z\in\ov{G}$ and $\xi\in\Om\setminus\{\infty\}$,
\begin{eqnarray*}
|\xi-z|&\le& |\xi-\xi_L|+|\xi_L-z|\\
&\le& 2|\xi_L-\xi|+|\xi-z|\ole|\xi-z|,
\end{eqnarray*}
i.e.,
\beq\label{2.3z}
|\xi-z|\asymp|\xi-\xi_L|+|\xi_L-z|.
\eeq
 Let for $0<\tau\le 1$ and $\z\in\ov{\Om^*_\de}\setminus\{\infty\}$,
$$
\tilde{\z}_{\de,\tau}:=\Psi_\de((1+\tau)\Phi_\de(\z)).
$$
Lemma \ref{lem2.3} with $F=\Phi_\de$ implies also that for $z\in \ov{\Om^*_\de}\setminus\{\infty\}$
and $\z\in\ov{\Om^*_\de}$ with $|\z-z|\le c_6|z-\tilde{z}_{\de,\tau}|$ we have
\beq\label{2.15}
c_7^{-1}|z-\tilde{z}_{\de,\tau}|\le |\z-\tilde{\z}_{\de,\tau}|\le
c_7 |z-\tilde{z}_{\de,\tau}|,
\eeq
where $c_7=c_7(K,c_6)>1$.

Furthermore, let $0<\de=\tau<1/2$. Since by \cite[p. 376, Lemma 2.2]{andbla}
and Lemma \ref{lem2.3},
\beq\label{2.6z}
|\Phi_\de(z)|-1\asymp\de,\quad z\in L ,
\eeq
 (\ref{2.2z}) and (\ref{2.4z}), written for $L_\de^*$ and $\Phi_\de$
instead of $L$ and $\Phi$, yield
\beq\label{2.1z}
d(z,L)\asymp|z-\tilde{z}_{\de,\de}|,\quad z\in L_\de^*.
\eeq
Moreover, we claim that
\beq\label{2.16}
|\tilde{z}_{\de,\de}-z|\asymp |\tilde{z}_{\de}-z|,\quad z\in L.
\eeq
Indeed, 
let $z^\bullet=z^\bullet(\de)\in L_\de^*$ satisfy $|z-z^\bullet|=d(z,L^*_\de)$.
Applying Lemma \ref{lem2.3} with $F=\Psi$ twice: first with the triplet
$\Phi(z),(1-\de)\Phi(z), \Phi(z^\bullet)$ and then 
with the triplet
$\Phi(z),(1-\de)\Phi(z), (1+\de)\Phi(z)$, we obtain
\beq\label{2.8z}
|z-\tilde{z}_\de|\asymp d(z,L_\de^*).
\eeq
Let $z'=z'(\de):=\Psi_\de(\Phi_\de(z)/|\Phi_\de(z)|)$
so that by (\ref{2.5z}), written for $\Om_\de^*$ instead of $\Om$,
we have
 $d(z,L_\de^*)\asymp|z-z'|$.
Since by (\ref{2.6z}) and Lemma \ref{lem2.3} with $F=\Psi_\de$ and
the triplet $\Phi_\de(z'),\Phi_\de(z), (1+\de)\Phi_\de(z)$
$$
|z'-z|\asymp |z'-\tilde{z'}_{\de,\de}|,
$$
according to (\ref{2.15})
\beq\label{2.9z}
|z-\tilde{z}_{\de,\de}|\asymp |z'-\tilde{z'}_{\de,\de}|\asymp
|z'-z|\asymp d(z,L_\de^*).
\eeq
Comparing (\ref{2.8z}) and (\ref{2.9z}) we obtain (\ref{2.16}).

To estimate the Christoffel function from above we use special
polynomials defined as follows.
For $\xi\in\Om\setminus\{\infty\}, z\in \ov{G}$, and $n\in\N$ with $n\ge 2$,
consider the {\it Dzjadyk kernel} $K_{0,1,1,n}(\xi,z)$ associated with $\ov{G}$
(see \cite[p. 429]{dzj} or \cite[p. 387]{andbla})  which is a polynomial in $z$ of degree
at most $4n$ with coefficients depending on $\xi$.
By virtue of \cite[p. 389, Theorem 2.4]{andbla}  we have
\beq\label{2.17}
\left|\frac{1}{\xi-z}-K_{0,1,1,n}(\xi,z)\right|\le \frac{c_8}{|\xi-z|}\left|
\frac{\tilde{\xi}_{1/n}-\xi}{\tilde{\xi}_{1/n}-z}\right|.
\eeq
For $\z\in L$, define  $\xi=\xi(\z,n):=\tilde{\z}_{c_9/n}, n>c_9$, where $c_9>1$ is chosen
as follows.
According to Lemma \ref{lem2.3} with $F=\Phi$ and the triplet
$\tilde{\xi}_{1/n},\xi,\z$ as well as
 (\ref{2.5z}), for $z\in \ov{G}$,
\beq\label{2.19}
 \left|
\frac{\tilde{\xi}_{1/n}-\xi}{\tilde{\xi}_{1/n}-z}\right|=
\left|
\frac{\tilde{\xi}_{1/n}-\xi}{\tilde{\xi}_{1/n}-\z}\right|
\left|
\frac{\tilde{\xi}_{1/n}-\z}{\tilde{\xi}_{1/n}-z}\right|\le \frac{c_{10}}{c_9^{1/K^2}}
<\frac{1}{2c_8}
\eeq
if $c_9:=1+(2c_8c_{10})^{K^2}$.

Since (\ref{2.17}) and (\ref{2.19}) imply for $z\in\ov{G}$,
$$
\left|\frac{1}{\xi-z}-K_{0,1,1,n}(\xi,z)\right|\le \frac{1}{2|\xi-z|},
$$
by (\ref{2.2z}), (\ref{2.4z}), and (\ref{2.3z}) we have
\beq\label{2.20}
|K_{0,1,1,n}(\xi,z)|\asymp |\xi-z|^{-1}\asymp
(|\z-z|+\rho_{1/n}(\z))^{-1}.
\eeq
For $\z\in L$ and any (fixed) $s\in\N$, consider polynomials (in $z$)
 of degree at most $4sn$ defined by
$$
q_{n,s,\z}(z):=(\rho_{1/n}(\z)K_{0,1,1,n}(\xi(\z,n),z))^s,\quad
Q_{n,s,\z}(z):=\frac{q_{n,s,\z}(\z,z)}{q_{n,s,\z}(\z,\z)}.
$$
Summarizing, we let
$$
p_{n,s,\z}
:= \left\{\begin{array}{ll}
   \displaystyle
 1
 &\mb{ if }n\le 8s,
\\[2ex]
Q_{\left\lfloor n/(4s)\right\rfloor,s,\z}&\mb{ if }n>8s,
\end{array}\right.
$$
where  $\left\lfloor x\right\rfloor$ denotes the {\it integer part} of a real number
$x$,
and use (\ref{2.4z}) and (\ref{2.20}) to obtain the following statement.
\begin{lem}\label{lem2.4}
For $n\in\N,\z\in L$ and  fixed $s\in\N$ there exists a polynomial
$p_{n,s,\z}\in\bP_n$ with the following properties:

(i) $p_{n,s,\z}(\z)=1$;

(ii) for $z\in \ov{G}$, 
$$
|p_{n,s,\z}(z)|\le c_{11}\left(\frac{\rho_{1/n}(\z)}{|\z-z|+\rho_{1/n}(\z)}\right)^s,
$$
where $c_{11}=c_{11}(G,s)$.
\end{lem}

\absatz{Proof of Theorems}

We start with a  modification of the classical Ahlfors result \cite{ahl1}.
As before, denote by $K\ge1$ a coefficient of quasiconformality of $L$.
\begin{lem}\label{lem3.0} (see \cite[pp. 25-26, Lemma 1.4 
and Corollary 1.3]{andbeldzj}).
There exists a quasiconformal reflection $y:\OC\to\OC$ with respect to $L$
satisfying the following properties:

(i) $y(G)=\Om,y(\Om)=G$,
$$
y(z)=z,\quad z\in L;
$$

(ii) $y$ has continuous partial derivatives of first order in
$\C\setminus(L\cup\{z_0\})$, where $z_0:=y(\infty)$;

(iii) for $\z_1,\z_2\in\ov{G}\setminus D_0$, where 
$D_0:=D(z_0,d(z_0,L)/2)$, the inequality
$$
c_1^{-1}|\z_1-\z_2|\le |y(\z_1)-y(\z_2)|\le c_1|\z_1-\z_2|
$$
holds with $c_1=c_1(K)>1$;

(iv) the inequalities
$$
|y_{\ov{\z}}(\z)|\le c_2|y(\z)|^2,\quad \z\in D_0,
$$
$$
|y_{\ov{\z}}(\z)|\le c_2,\quad \z\in G\setminus D_0,
$$
hold with $c_2=c_2(K)$.
\end{lem}
Next, we claim that 
\beq\label{3.2}
|z-\z|+d(z,L)\ole |z-y(\z)|,\quad z,\z\in G\setminus D_0.
\eeq
Indeed, in the nontrivial case where $d(z,L)<|z-\z|$, we introduce a point
$z'\in L$ such that $d(z,L)=|z-z'|$ and use Lemma \ref{lem3.0} to obtain
\begin{eqnarray*}
|z-\z|+d(z,L)&<&2|z-\z|\le 2|z-z'|+ 2|z'-\z|\ole
|z-z'|+|z'-y(\z)|\\
&\le& 2|z-z'|+|z-y(\z)|\le 3|z-y(\z)|
\end{eqnarray*}
which proves (\ref{3.2}).

For the weight function $h$ defined by (\ref{1.3}) and $1\le p<\infty$, denote by 
$A_p(h,G)$ the space of  functions $f$ analytic in $G$ and satisfying
$$
||f||_{A_p(h,G)}^p:=\int_G|f|^phdm<\infty.
$$
Note that polynomials are in $A_p(h,G)$.
\begin{lem}\label{lem3.1}
For $p_n\in\bP_n,n\in\N,1\le p<\infty$, and $z\in G\setminus D_0$,
we have
\beq\label{3.3}
|p_n'(z)|\le c_3d(z,L)^{-1-2/p}\prod_{j=1}^m|z-z_j|^{-\al_j/p}
||p_n||_{A_p(h,G)},
\eeq
where $c_3=c_3(G,h,p)$.
\end{lem}
{\bf Proof.} Consider an analytic in $G$ function
$$
H_p(z):=\prod_{j=1}^m(z-z_j)^{\al_j/p}.
$$
Since $\int_G|p_nH_p|dm<\infty$, we can use the  Belyi integral
formula (see \cite{bel} or \cite[p. 110, Theorem 4.4]{andbeldzj})
to obtain for $z\in G\setminus D_0$
\begin{eqnarray}
p_n'(z)H_p(z)&=&-p_n(z)H_p'(z)+(p_n(z)H_p(z))'\nonumber\\
&=&-p_n(z)H_p'(z)-\frac{2}{\pi}\int_G\frac{p_n(\z)H_p(\z)}{(y(\z)-z)^3}
y_{\ov{\z}}(\z)dm(\z)\nonumber\\
\label{3.4}
&=:&-A(z)-B(z).
\end{eqnarray}
According to the mean-value property for a subharmonic function $|p_n|^p$
(see \cite[p. 46, Theorem 2.6.8(b)]{ran}), letting $d:=d(z,L)$ we have
$$
|p_n(z)|^p\le\frac{4}{\pi d^2}\int_{D(z,d/2)}|p_n|^pdm
\ole
d^{-2}\prod_{j=1}^m|z-z_j|^{-\al_j}||p_n||_{A_p(h,G)}^p
$$
which implies
\begin{eqnarray}
|A(z)|&=&|p_n(z)||H_p(z)|\left|\frac{H_p'(z)}{H_p(z)}\right|
\nonumber\\
\label{3.5}
&\le&
|p_n(z)||H_p(z)|\sum_{j=1}^m\frac{|\al_j|}{p|z-z_j|}\ole
d^{-1-2/p}
||p_n||_{A_p(h,G)}.
\end{eqnarray}
To estimate $|B(z)|$ we consider two particular cases. 

If $p=1$, then by Lemma \ref{lem3.0} and (\ref{3.2}),
\begin{eqnarray}
|B(z)|&\ole&
\int_{D_0}|p_n|hdm+d^{-3}\int_{G\setminus D_0}|p_n|h
dm\nonumber\\
\label{3.6}
&\ole& d^{-3}||p_n||_{A_1(h,G)}.
\end{eqnarray}
If $p>1$,  then H\"older's inequality with $q:=p/(p-1)$ yields
$$
|B(z)|\ole ||p_n||_{A_p(h,G)}\left(\int_G\frac{|y_{\ov{\z}}(\z)|^qdm(\z)}{|y(\z)-z|^{3q}}
\right)^{1/q}=:||p_n||_{A_p(h,G)}C(z)^{1/q}.
$$
According to Lemma \ref{lem2.2} with $\be=3q,\de=d$ and $\al_j=0$ as well as Lemma \ref{lem3.0}
and (\ref{3.2}) we obtain
$$
C(z)\ole\int_{D_0}dm+
\int_{G\setminus D_0}\frac{dm(\z)}{(|\z-z|+d)^{3q}}\ole
d^{2-3q}
$$
which implies
\beq\label{3.7}
|B(z)|\ole d^{-1-2/p}||p_n||_{A_p(h,G)}.
\eeq
Comparing (\ref{3.4})-(\ref{3.7}) we have (\ref{3.3}).

\hfill$\Box$

\begin{lem}\label{lem3.2}
There exists $\ve=\ve(G)$ such that for $p_n\in \bP_n, n\in\N,
1\le p<\infty, z\in L$, and $\z\in D(z,\ve \rho_{1/n}(z))$
the inequality
\beq\label{3.13}
|p_n'(\z)|\le c_4\rho_{1/n}(z)^{-1-2/p}\prod_{j=1}^m(\rho_{1/n}(z)
+|z-z_j|)^{-\al_j/p}||p_n||_{A_p(h,G)}
\eeq
holds with $c_4=c_4(\ve, G,h,p)$.
\end{lem}
{\bf Proof.}  Without loss of generality we assume that
 $n>2$ and, in addition to the points $z_j\in L, j=1,\ldots,m$, we 
introduce points 
$$
z_{j,n}:=\Psi\left(\left(1-\frac{2}{n}\right)\Phi(z_j)\right)
$$
which, according to Lemma \ref{lem2.3}
with $F=\Phi$ and the triplet of points $z,z_j,z_{j,n}$
 satisfy
\beq\label{3.8}
|z-z_j|\asymp|z-z_{j,n}|,\quad z\in L^*_{1/n}.
\eeq
Let $1\le k\le n$ be the degree of $p_n$ and let
$$
\z_n:=\tilde{\z}_{1/n,1/n}=\Psi_{1/n}\left(\left(1+\frac{1}{n}\right)
\Phi_{1/n}(\z)\right),\quad \z\in\Om^*_{1/n}.
$$
Consider subharmonic in $\Om^*_{1/n}$ function
\begin{eqnarray*}
f(\z)=f_{n,p}(\z)&:=&
\ln|p_n'(\z)|+\left(1+\frac{2}{p}\right)\ln|\z-\z_n|
+\frac{1}{p}\sum_{j=1}^m\al_j\ln|\z-z_{j,n}|\\
&&-\left( k+\frac{2}{p}+\frac{1}{p}\sum_{j=1}^m\al_j\right)
\ln|\Phi_{1/n}(\z)|
\end{eqnarray*}
(which is harmonic and bounded in a neighborhood of infinity).

Since by virtue of (\ref{2.1z}), Lemma \ref{lem3.1}, and (\ref{3.8}), 
$$
f(\z)\le c_5+\ln ||p_n||_{A_p(h,G)},\quad \z\in L^*_{1/n},
$$
by the maximum principle (see \cite[p. 29]{ran}) the same inequality
holds for $\z\in D(z,\ve \rho_{1/n}(z))$, where $z\in L$ and $\ve$ is
chosen such that
$$
\ve\rho_{1/n}(z)\le\frac{1}{2}d(z,L^*_{1/n}).
$$
The existence of such a constant $\ve$ is guaranteed by Lemma \ref{lem2.3}
with $F=\Phi$ and (\ref{2.2z}).

Applying (\ref{2.nn}) and  Lemma \ref{lem2.3} 
with $F=\Phi$ and the triplet $z,z_{j,n}, \tilde{(z_j)}_{1/n}$
we obtain
for $\z\in D(z,\ve \rho_{1/n}(z))$
$$
|\z-z_{j,n}|\asymp|z-z_{j,n}|
\asymp|z-\tilde{(z_j)}_{1/n}|
\asymp\rho_{1/n}(z)+|z-z_j|.
$$
Therefore, according to (\ref{2.2z}), 
 (\ref{2.15}), (\ref{2.6z}),
and (\ref{2.16}) we further have
\begin{eqnarray*}
\ln|p_n'(\z)|&=& f(\z)-
\left(1+\frac{2}{p}\right)\ln|\z-\z_n|
-\frac{1}{p}\sum_{j=1}^m\al_j\ln|\z-z_{j,n}|\\
&&+\left( k+\frac{2}{p}+\frac{1}{p}\sum_{j=1}^m\al_j\right)
\ln|\Phi_{1/n}(\z)|\\
&\le& c_6 +
\ln ||p_n||_{A_p(h,G)} -\left(1+\frac{2}{p}\right)\ln\rho_{1/n}(z)\\
&&
-\frac{1}{p}\sum_{j=1}^m\al_j\ln(\rho_{1/n}(z)+|z-z_j|),
\end{eqnarray*}
which
yields (\ref{3.13}).

\hfill$\Box$

{\bf Proof of Theorem \ref{th1}.} 
Let $s\in\N$ be a fixed number with $s>2+\sum_{j=1}^m|\al_j|$.
Consider polynomial $p_n=p_{n,s,z}$ from Lemma \ref{lem2.4}.
On account of (\ref{1.1}), Lemma \ref{lem2.2} with $\be:=sp$,  $\de:=\rho_{1/n}(z)$
and Lemma \ref{lem2.4}, we have
\begin{eqnarray*}
\la_n(\nu,p,z)&\le&\int_{G}|p_n|^phdm\\
&\ole&
\rho_{1/n}(z)^{sp}\int_G(|z-\z|+\rho_{1/n}(z))^{-sp}\prod_{j=1}^m
|\z-z_j|^{\al_j}dm(\z)\\
&\ole&
\rho_{1/n}(z)^{2}\prod_{j=1}^m(|z-z_j|+\rho_{1/n}(z))^{\al_j},
\end{eqnarray*}
which proves the right-hand side of (\ref{1.6}).

In order to prove the left-hand side of (\ref{1.6}), it is sufficient
to show that for any $z\in L$ and $p_n\in\bP_n$ with $p_n(z)=1$,
we have
\beq\label{3.21}
\int_G|p_n|^phdm\oge
\rho_{1/n}(z)^{2}\prod_{j=1}^m(|z-z_j|+\rho_{1/n}(z))^{\al_j}.
\eeq
In the nontrivial  case where
$$
\int_G|p_n|^phdm\le
\rho_{1/n}(z)^{2}\prod_{j=1}^m(|z-z_j|+\rho_{1/n}(z))^{\al_j}
$$
 Lemma \ref{lem3.2} implies 
$$
|p_n'(\z)|\le \frac{c_4}{\rho_{1/n}(z)}, \quad
\z\in D(z,\ve\rho_{1/n}(z)).
$$
Moreover, for the same $\z$
$$
|p_n(\z)-1|\le\int_{[z,\z]}|p_n'(\xi)||d\xi|
\le c_4\frac{|z-\z|}{\rho_{1/n}(z)}
$$
and if 
$$
|z-\z|\le \ve_1\rho_{1/n}(z),\quad \ve_1:=\min\left(\ve,\frac{1}{2c_4}\right),
$$
we obtain $|p_n(\z)|\ge1/2$.

Hence, according to (\ref{2.7}) with $\de=\ve_2\rho_{1/n}(z)$,
where $\ve_2=\ve_2(G)<\ve_1$ is chosen so that $\de<\min_{j\neq k}|z_j-z_k|/4$,
we have
\begin{eqnarray*}
\int_G|p_n|^phdm&\ge&
\int_{G\cap D(z,\ve_2\rho_{1/n}(z))}|p_n|^phdm\\
&\oge&
\rho_{1/n}(z)^{2}\prod_{j=1}^m(|z-z_j|+\rho_{1/n}(z))^{\al_j}
\end{eqnarray*}
which proves (\ref{3.21}).

\hfill$\Box$

{\bf Proof of Theorem \ref{th2}}.
By Lemma \ref{lem2.4} with the quasidisk 
$$
G=\{z=x+iy:0<x<1,|y|< e^{-1}\}
$$
for any fixed $k\in\N$ and any
$n\in\N$ there
exists a polynomial $p_n=p_{n,4k,0}\in\bP_n$
such that $p_n(0)=1$ and
$$
|p_n(\z)|
\ole \left\{\begin{array}{ll}
   \displaystyle
 1
 &\mb{ if }\z\in G, |\z|\le 1/n,
\\[2ex]
(n|\z|)^{-4k}&\mb{ if }\z\in G, |\z|>1/n.
\end{array}\right.
$$
Since $G^*\subset G$ and by the 
L\"owner inequality \cite{loe}
  (see also \cite[p. 359, Corollary 2.5]{andbla}),
	$$
  \rho_{1/n}(0)
  \ge \frac{\mb{diam }G^*}{8n^2},
 $$
by virtue of ({1.1}) we obtain
\begin{eqnarray*}
\la_n(m^*,p,0)&\le&
\int_{G^*}|p_n|^pdm^*\le
\int_{G^*\cap D(0,1/n)}|p_n|^pdm^*
+\int_{G^*\setminus D(0,1/n)}|p_n|^pdm^*\\
&\ole&\frac{1}{n}e^{-n}+n^{-4kp}\int_{1/(2n)}^1 x^{-4kp}
e^{-1/x}dx\ole n^{-4k}\ole\rho_{1/n}(0)^{2k}
\end{eqnarray*}
from which (\ref{1.9}) follows.

\hfill$\Box$

{\bf Acknowledgements.}
Part of this work was done during the Fall of 2016 semester, while the author visited
the Katholische Universit\"at Eichst\"att-Ingolstadt and
the Julius Maximilian University of W\"urzburg.
The author is  also grateful to   M.
 Nesterenko
 for his helpful comments.

V. V. Andrievskii

 Department of Mathematical Sciences

 Kent State University

 Kent, OH 44242

 USA

e-mail: andriyev@math.kent.edu

tel: 330-672-9029

\end{document}